\documentclass{amsart}
\usepackage{amssymb}
\usepackage{amsfonts}

\newtheorem{theorem}{Theorem}
\theoremstyle{plain}
\newtheorem{acknowledgement}{Acknowledgement}

\newtheorem{corollary}{Corollary}

\newtheorem{remark}{Remark}

\numberwithin{equation}{section}
\input{tcilatex}

\begin{document}
\title[genearlizations of Precupanu's Inequality ]{Genearlizations of
Precupanu's Inequality for Orthornormal Families of Vectors in Inner Product
Spaces }
\author{S.S. Dragomir}
\address{School of Computer Science and Mathematics\\
Victoria University of Technology\\
PO Box 14428, MCMC 8001\\
VIC, Australia.}
\email{sever@matilda.vu.edu.au}
\urladdr{http://rgmia.vu.edu.au/SSDragomirWeb.html}
\date{November 10, 2004.}

\begin{abstract}
Some genearlizations of Precupanu's inequality for orthornormal families of
vectors in real or complex inner product spaces and applications related to
Buzano's, Richard's and Kurepa's results are given.
\end{abstract}

\maketitle

\section{\label{s1}Introduction}

In 1976, T. Precupanu \cite{P} obtained the following result related to the
Schwarz inequality in a real inner product space $\left( H;\left\langle
\cdot ,\cdot \right\rangle \right) :$

\begin{theorem}
\label{t1.1}For any $a\in H,$ $x,y\in H\backslash \left\{ 0\right\} ,$ we
have the inequality:%
\begin{align}
\frac{-\left\Vert a\right\Vert \left\Vert b\right\Vert +\left\langle
a,b\right\rangle }{2}& \leq \frac{\left\langle x,a\right\rangle \left\langle
x,b\right\rangle }{\left\Vert x\right\Vert ^{2}}+\frac{\left\langle
y,a\right\rangle \left\langle y,b\right\rangle }{\left\Vert y\right\Vert ^{2}%
}-2\cdot \frac{\left\langle x,a\right\rangle \left\langle y,b\right\rangle
\left\langle x,y\right\rangle }{\left\Vert x\right\Vert ^{2}\left\Vert
y\right\Vert ^{2}}  \label{1.1} \\
& \leq \frac{\left\Vert a\right\Vert \left\Vert b\right\Vert +\left\langle
a,b\right\rangle }{2}.  \notag
\end{align}%
In the right-hand side or in the left-hand side of (\ref{1.1}) we have
equality if and only if there are $\lambda ,\mu \in \mathbb{R}$ such that%
\begin{equation}
\lambda \frac{\left\langle x,a\right\rangle }{\left\Vert x\right\Vert ^{2}}%
\cdot x+\mu \frac{\left\langle y,b\right\rangle }{\left\Vert y\right\Vert
^{2}}\cdot y=\frac{1}{2}\left( \lambda a+\mu b\right) .  \label{1.2}
\end{equation}
\end{theorem}

Note for instance that \cite{P}, if $y\perp b,$ i.e., $\left\langle
y,b\right\rangle =0,$ then by (\ref{1.1}) one may deduce:%
\begin{equation}
\frac{-\left\Vert a\right\Vert \left\Vert b\right\Vert +\left\langle
a,b\right\rangle }{2}\left\Vert x\right\Vert ^{2}\leq \left\langle
x,a\right\rangle \left\langle x,b\right\rangle \leq \frac{\left\Vert
a\right\Vert \left\Vert b\right\Vert +\left\langle a,b\right\rangle }{2}%
\left\Vert x\right\Vert ^{2}  \label{1.3}
\end{equation}%
for any $a,b,x\in H,$ an inequality that has been obtained previously by U.
Richard \cite{R}. The case of equality in the right-hand side or in the
left-hand side of (\ref{1.3}) holds if and only if there are $\lambda ,\mu
\in \mathbb{R}$ with%
\begin{equation}
2\lambda \left\langle x,a\right\rangle x=\left( \lambda a+\mu b\right)
\left\Vert x\right\Vert ^{2}.  \label{1.4}
\end{equation}%
For $a=b,$ we may obtain from (\ref{1.1}) the following inequality \cite{P}%
\begin{equation}
0\leq \frac{\left\langle x,a\right\rangle ^{2}}{\left\Vert x\right\Vert ^{2}}%
+\frac{\left\langle y,a\right\rangle ^{2}}{\left\Vert y\right\Vert ^{2}}%
-2\cdot \frac{\left\langle x,a\right\rangle \left\langle y,a\right\rangle
\left\langle x,y\right\rangle }{\left\Vert x\right\Vert ^{2}\left\Vert
y\right\Vert ^{2}}\leq \left\Vert a\right\Vert ^{2}.  \label{1.5}
\end{equation}%
This inequality implies \cite{P}:%
\begin{equation}
\frac{\left\langle x,y\right\rangle }{\left\Vert x\right\Vert \left\Vert
y\right\Vert }\geq \frac{1}{2}\left[ \frac{\left\langle x,a\right\rangle }{%
\left\Vert x\right\Vert \left\Vert a\right\Vert }+\frac{\left\langle
y,a\right\rangle }{\left\Vert y\right\Vert \left\Vert a\right\Vert }\right]
^{2}-\frac{3}{2}.  \label{1.6}
\end{equation}

In \cite{M}, M.H. Moore pointed out the following reverse of the Schwarz
inequality%
\begin{equation}
\left\vert \left\langle y,z\right\rangle \right\vert \leq \left\Vert
y\right\Vert \left\Vert z\right\Vert ,\qquad y,z\in H,  \label{1.7}
\end{equation}%
where some information about a third vector $x$ is known:

\begin{theorem}
\label{t1.2}Let $\left( H;\left\langle \cdot ,\cdot \right\rangle \right) $
be an inner product space over the real field $\mathbb{R}$ and $x,y,z\in H$
such that:%
\begin{equation}
\left\vert \left\langle x,y\right\rangle \right\vert \geq \left(
1-\varepsilon \right) \left\Vert x\right\Vert \left\Vert y\right\Vert
,\qquad \left\vert \left\langle x,z\right\rangle \right\vert \geq \left(
1-\varepsilon \right) \left\Vert x\right\Vert \left\Vert z\right\Vert ,
\label{1.8}
\end{equation}%
where $\varepsilon $ is a positive real number, reasonably small. Then%
\begin{equation}
\left\vert \left\langle y,z\right\rangle \right\vert \geq \max \left\{
1-\varepsilon -\sqrt{2\varepsilon },1-4\varepsilon ,0\right\} \left\Vert
y\right\Vert \left\Vert z\right\Vert .  \label{1.9}
\end{equation}
\end{theorem}

Utilising Richard's inequality (\ref{1.3}) written in the following
equivalent form:%
\begin{equation}
2\cdot \frac{\left\langle x,a\right\rangle \left\langle x,b\right\rangle }{%
\left\Vert x\right\Vert ^{2}}-\left\Vert a\right\Vert \left\Vert
b\right\Vert \leq \left\langle a,b\right\rangle \leq 2\cdot \frac{%
\left\langle x,a\right\rangle \left\langle x,b\right\rangle }{\left\Vert
x\right\Vert ^{2}}+\left\Vert a\right\Vert \left\Vert b\right\Vert
\label{1.10}
\end{equation}%
for any $a,b\in H$ and $a\in H\backslash \left\{ 0\right\} ,$ Precupanu has
obtained the following Moore's type result:

\begin{theorem}
\label{t1.3}Let $\left( H;\left\langle \cdot ,\cdot \right\rangle \right) $
be a real inner product space. If $a,b,x\in H$ and $0<\varepsilon
_{1}<\varepsilon _{2}$ are such that:%
\begin{align}
\varepsilon _{1}\left\Vert x\right\Vert \left\Vert a\right\Vert & \leq
\left\langle x,a\right\rangle \leq \varepsilon _{2}\left\Vert x\right\Vert
\left\Vert a\right\Vert ,  \label{1.11} \\
\varepsilon _{1}\left\Vert x\right\Vert \left\Vert b\right\Vert & \leq
\left\langle x,b\right\rangle \leq \varepsilon _{2}\left\Vert x\right\Vert
\left\Vert b\right\Vert ,  \notag
\end{align}%
then%
\begin{equation}
\left( 2\varepsilon _{1}^{2}-1\right) \left\Vert a\right\Vert \left\Vert
b\right\Vert \leq \left\langle a,b\right\rangle \leq \left( 2\varepsilon
_{1}^{2}+1\right) \left\Vert a\right\Vert \left\Vert b\right\Vert .
\label{1.12}
\end{equation}
\end{theorem}

Remark that the right inequality is always satisfied, since by Schwarz's
inequality, we have $\left\langle a,b\right\rangle \leq \left\Vert
a\right\Vert \left\Vert b\right\Vert $. The left inequality may be useful
when one assumes that $\varepsilon _{1}\in (0,1].$ In that case, from (\ref%
{1.12}), we obtain%
\begin{equation}
-\left\Vert a\right\Vert \left\Vert b\right\Vert \leq \left( 2\varepsilon
_{1}^{2}-1\right) \left\Vert a\right\Vert \left\Vert b\right\Vert \leq
\left\langle a,b\right\rangle  \label{1.13}
\end{equation}%
provided $\varepsilon _{1}\left\Vert x\right\Vert \left\Vert a\right\Vert
\leq \left\langle x,a\right\rangle $ and $\varepsilon _{1}\left\Vert
x\right\Vert \left\Vert b\right\Vert \leq \left\langle x,b\right\rangle ,$
which is a refinement of Schwarz's inequality%
\begin{equation*}
-\left\Vert a\right\Vert \left\Vert b\right\Vert \leq \left\langle
a,b\right\rangle .
\end{equation*}

In the complex case, apparently independent of Richard, M.L. Buzano obtained
in \cite{B} the following inequality%
\begin{equation}
\left\vert \left\langle x,a\right\rangle \left\langle x,b\right\rangle
\right\vert \leq \frac{\left\Vert a\right\Vert \left\Vert b\right\Vert
+\left\vert \left\langle a,b\right\rangle \right\vert }{2}\cdot \left\Vert
x\right\Vert ^{2},  \label{1.14}
\end{equation}%
provided $x,a,b$ are vectors in the complex inner product space $\left(
H;\left\langle \cdot ,\cdot \right\rangle \right) .$

In the same paper \cite{P}, Precupanu, without mentioning Buzano's name in
relation to the inequality (\ref{1.14}), observed that, on utilising (\ref%
{1.14}), one may obtain the following result of Moore type:

\begin{theorem}
\label{t1.4}Let $\left( H;\left\langle \cdot ,\cdot \right\rangle \right) $
be a (real or) complex inner product space. If $x,a,b\in H$ are such that%
\begin{equation}
\left\vert \left\langle x,a\right\rangle \right\vert \geq \left(
1-\varepsilon \right) \left\Vert x\right\Vert \left\Vert a\right\Vert
,\qquad \left\vert \left\langle x,b\right\rangle \right\vert \geq \left(
1-\varepsilon \right) \left\Vert x\right\Vert \left\Vert b\right\Vert ,
\label{1.15}
\end{equation}%
then%
\begin{equation}
\left\vert \left\langle a,b\right\rangle \right\vert \geq \left(
1-4\varepsilon +2\varepsilon ^{2}\right) \left\Vert a\right\Vert \left\Vert
b\right\Vert .  \label{1.16}
\end{equation}
\end{theorem}

Note that the above theorem is useful when, for $\varepsilon \in (0,1],$ the
quantity $1-4\varepsilon +2\varepsilon ^{2}>0,$ i.e., $\varepsilon \in
\left( 0,1-\frac{\sqrt{2}}{2}\right] .$

\begin{remark}
When the space is real, the inequality (\ref{1.16}) provides a better lower
bound for $\left\vert \left\langle a,b\right\rangle \right\vert $ than the
second bound in Moore's result (\ref{1.9}). However, it is not known if the
first bound in (\ref{1.9}) remains valid for the case of complex spaces.
From Moore's original proof, apparently, the fact that the space $\left(
H;\left\langle \cdot ,\cdot \right\rangle \right) $ is real plays an
essential role.
\end{remark}

Before we point out some new results for orthonormal families of vectors in
real or complex inner product spaces, we state the following result that
complements the Moore type results outlined above for real spaces:

\begin{theorem}
\label{t1.5}Let $\left( H;\left\langle \cdot ,\cdot \right\rangle \right) $
be a real inner product space and $a,b,x,y\in H\backslash \left\{ 0\right\}
. $

\begin{enumerate}
\item[(i)] If there exist $\delta _{1},\delta _{2}\in (0,1]$ such that%
\begin{equation*}
\frac{\left\langle x,a\right\rangle }{\left\Vert x\right\Vert \left\Vert
a\right\Vert }\geq \delta _{1},\qquad \frac{\left\langle y,a\right\rangle }{%
\left\Vert y\right\Vert \left\Vert a\right\Vert }\geq \delta _{2}
\end{equation*}%
and $\delta _{1}+\delta _{2}\geq 1,$ then%
\begin{equation}
\frac{\left\langle x,y\right\rangle }{\left\Vert x\right\Vert \left\Vert
y\right\Vert }\geq \frac{1}{2}\left( \delta _{1}+\delta _{2}\right) ^{2}-%
\frac{3}{2}\qquad \left( \geq -1\right) .  \label{1.17}
\end{equation}

\item[(ii)] If there exist $\mu _{1}\left( \mu _{2}\right) \in \mathbb{R}$
such that%
\begin{equation*}
\mu _{1}\left\Vert a\right\Vert \left\Vert b\right\Vert \leq \frac{%
\left\langle x,a\right\rangle \left\langle x,b\right\rangle }{\left\Vert
x\right\Vert ^{2}}\left( \leq \mu _{2}\left\Vert a\right\Vert \left\Vert
b\right\Vert \right)
\end{equation*}%
and $1\geq \mu _{1}\geq 0$ $\left( -1\leq \mu _{2}\leq 0\right) ,$ then%
\begin{equation}
\left[ -1\leq \right] 2\mu _{1}-1\leq \frac{\left\langle a,b\right\rangle }{%
\left\Vert a\right\Vert \left\Vert b\right\Vert }\left( \leq 2\mu _{2}+1%
\left[ \leq 1\right] \right) .  \label{1.18}
\end{equation}
\end{enumerate}
\end{theorem}

The proof is obvious by the inequalities (\ref{1.6}) and (\ref{1.10}). We
omit the details.

\section{Inequalities for orthonormal Families}

We recall that the finite family $\left\{ e_{i}\right\} _{i\in I}$ is
orthonormal in $\left( H;\left\langle \cdot ,\cdot \right\rangle \right) ,$
a real or complex inner product space, if 
\begin{equation*}
\left\langle e_{i},e_{j}\right\rangle =\left\{ 
\begin{array}{ll}
0 & \text{if \ }i\neq j \\ 
&  \\ 
1 & \text{if \ }i=j%
\end{array}%
\right.
\end{equation*}%
where $i,j\in I.$

The following result may be stated.

\begin{theorem}
\label{t2.1}Let $\left\{ e_{i}\right\} _{i\in I}$ and $\left\{ f_{j}\right\}
_{j\in J}$ be two finite families of orthonormal vectors in $\left(
H;\left\langle \cdot ,\cdot \right\rangle \right) .$ For any $x,y\in
H\backslash \left\{ 0\right\} $ one has the inequality%
\begin{multline}
\left\vert \sum_{i\in I}\left\langle x,e_{i}\right\rangle \left\langle
e_{i},y\right\rangle +\sum_{j\in J}\left\langle x,f_{j}\right\rangle
\left\langle f_{j},y\right\rangle \right.  \label{2.1} \\
-\left. 2\sum_{i\in I,j\in J}\left\langle x,e_{i}\right\rangle \left\langle
f_{j},y\right\rangle \left\langle e_{i},f_{j}\right\rangle -\frac{1}{2}%
\left\langle x,y\right\rangle \right\vert \leq \frac{1}{2}\left\Vert
x\right\Vert \left\Vert y\right\Vert .
\end{multline}%
The case of equality holds in (\ref{2.1}) if and only if there exists a $%
\lambda \in \mathbb{K}$ such that%
\begin{equation}
x-\lambda y=2\left( \sum_{i\in I}\left\langle x,e_{i}\right\rangle
e_{i}-\lambda \sum_{j\in J}\left\langle y,f_{j}\right\rangle f_{j}\right) .
\label{2.2}
\end{equation}
\end{theorem}

\begin{proof}
We know that, if $u,v\in H,$ $v\neq 0,$ then%
\begin{equation}
\left\Vert u-\frac{\left\langle u,v\right\rangle }{\left\Vert v\right\Vert
^{2}}\cdot v\right\Vert ^{2}=\frac{\left\Vert u\right\Vert ^{2}\left\Vert
v\right\Vert ^{2}-\left\vert \left\langle u,v\right\rangle \right\vert ^{2}}{%
\left\Vert v\right\Vert ^{2}}  \label{2.3}
\end{equation}%
showing that, in Schwarz's inequality%
\begin{equation}
\left\vert \left\langle u,v\right\rangle \right\vert ^{2}\leq \left\Vert
u\right\Vert ^{2}\left\Vert v\right\Vert ^{2},  \label{2.4}
\end{equation}%
the case of equality, for $v\neq 0,$ holds if and only if%
\begin{equation}
u=\frac{\left\langle u,v\right\rangle }{\left\Vert v\right\Vert ^{2}}\cdot
v,\   \label{2.5}
\end{equation}%
i.e. there exists a $\lambda \in \mathbb{R}$ \ such that $u=\lambda v.$

Now, let $u:=2\sum_{i\in I}\left\langle x,e_{i}\right\rangle e_{i}-x$ and $%
v:=2\sum_{j\in J}\left\langle y,f_{j}\right\rangle f_{j}-y.$

Observe that%
\begin{align*}
\left\Vert u\right\Vert ^{2}& =\left\Vert 2\sum_{i\in I}\left\langle
x,e_{i}\right\rangle e_{i}\right\Vert ^{2}-4\func{Re}\left\langle \sum_{i\in
I}\left\langle x,e_{i}\right\rangle e_{i},x\right\rangle +\left\Vert
x\right\Vert ^{2} \\
& =4\sum_{i\in I}\left\vert \left\langle x,e_{i}\right\rangle \right\vert
^{2}-4\sum_{i\in I}\left\vert \left\langle x,e_{i}\right\rangle \right\vert
^{2}+\left\Vert x\right\Vert ^{2}=\left\Vert x\right\Vert ^{2},
\end{align*}%
and, similarly%
\begin{equation*}
\left\Vert v\right\Vert ^{2}=\left\Vert y\right\Vert ^{2}.
\end{equation*}%
Also,%
\begin{multline*}
\left\langle u,v\right\rangle =4\sum_{i\in I,j\in J}\left\langle
x,e_{i}\right\rangle \left\langle f_{j},y\right\rangle \left\langle
e_{i},f_{j}\right\rangle +\left\langle x,y\right\rangle \\
-2\sum_{i\in I}\left\langle x,e_{i}\right\rangle \left\langle
e_{i},y\right\rangle -2\sum_{j\in J}\left\langle x,f_{j}\right\rangle
\left\langle f_{j},y\right\rangle .
\end{multline*}%
Therefore, by Schwarz's inequality (\ref{2.4}) we deduce the desired
inequality (\ref{2.1}). By (\ref{2.5}), the case of equality holds in (\ref%
{2.1}) if and only if there exists a $\lambda \in \mathbb{K}$ such that%
\begin{equation*}
2\sum_{i\in I}\left\langle x,e_{i}\right\rangle e_{i}-x=\lambda \left(
2\sum_{j\in J}\left\langle y,f_{j}\right\rangle f_{j}-y\right) ,
\end{equation*}%
which is equivalent to (\ref{2.2}).
\end{proof}

\begin{remark}
If in (\ref{2.2}) we choose $x=y,$ then we get the inequality:%
\begin{multline}
\left\vert \sum_{i\in I}\left\vert \left\langle x,e_{i}\right\rangle
\right\vert ^{2}+\sum_{j\in J}\left\vert \left\langle x,f_{j}\right\rangle
\right\vert ^{2}-2\sum_{i\in I,j\in J}\left\langle x,e_{i}\right\rangle
\left\langle f_{j},x\right\rangle \left\langle e_{i},f_{j}\right\rangle -%
\frac{1}{2}\left\Vert x\right\Vert ^{2}\right\vert  \label{2.6} \\
\leq \frac{1}{2}\left\Vert x\right\Vert ^{2}
\end{multline}%
for any $x\in H.$

If in the above theorem we assume that $I=J$ and $f_{i}=e_{i},$ $i\in I,$
then we get from (\ref{2.1}) the Schwarz inequality $\left\vert \left\langle
x,y\right\rangle \right\vert \leq \left\Vert x\right\Vert \left\Vert
y\right\Vert .$

If $I\cap J=\varnothing ,$ $I\cup J=K,$ $g_{k}=e_{k},$ $k\in I,$ $%
g_{k}=f_{k},$ $k\in J$ and $\left\{ g_{k}\right\} _{k\in K}$ is orthonormal,
then from (\ref{2.1}) we get:%
\begin{equation}
\left\vert \sum_{k\in K}\left\langle x,g_{k}\right\rangle \left\langle
g_{k},y\right\rangle -\frac{1}{2}\left\langle x,y\right\rangle \right\vert
\leq \frac{1}{2}\left\Vert x\right\Vert \left\Vert y\right\Vert ,\qquad
x,y\in H  \label{2.7}
\end{equation}%
which has been obtained earlier by the author in \cite{SSD}.
\end{remark}

If $I$ and $J$ reduce to one element, namely $e_{1}=\frac{e}{\left\Vert
e\right\Vert },$ $f_{1}=\frac{f}{\left\Vert f\right\Vert }$ with $e,f\neq 0,$
then from (\ref{2.1}) we get%
\begin{multline}
\left\vert \frac{\left\langle x,e\right\rangle \left\langle e,y\right\rangle 
}{\left\Vert e\right\Vert ^{2}}+\frac{\left\langle x,f\right\rangle
\left\langle f,y\right\rangle }{\left\Vert f\right\Vert ^{2}}-2\cdot \frac{%
\left\langle x,e\right\rangle \left\langle f,y\right\rangle \left\langle
e,f\right\rangle }{\left\Vert e\right\Vert ^{2}\left\Vert f\right\Vert ^{2}}-%
\frac{1}{2}\left\langle x,y\right\rangle \right\vert  \label{2.8} \\
\leq \frac{1}{2}\left\Vert x\right\Vert \left\Vert y\right\Vert ,\qquad
x,y\in H
\end{multline}%
which is the corresponding complex version of Precupanu's inequality (\ref%
{1.1}).

If in (\ref{2.8}) we assume that $x=y,$ then we get%
\begin{equation}
\left\vert \frac{\left\vert \left\langle x,e\right\rangle \right\vert ^{2}}{%
\left\Vert e\right\Vert ^{2}}+\frac{\left\vert \left\langle x,f\right\rangle
\right\vert ^{2}}{\left\Vert f\right\Vert ^{2}}-2\cdot \frac{\left\langle
x,e\right\rangle \left\langle f,e\right\rangle \left\langle e,f\right\rangle 
}{\left\Vert e\right\Vert ^{2}\left\Vert f\right\Vert ^{2}}-\frac{1}{2}%
\left\Vert x\right\Vert ^{2}\right\vert \leq \frac{1}{2}\left\Vert
x\right\Vert ^{2}.  \label{2.9}
\end{equation}

The following corollary may be stated:

\begin{corollary}
\label{c2.3}With the assumptions of Theorem \ref{t2.1}, we have:%
\begin{align}
& \left\vert \sum_{i\in I}\left\langle x,e_{i}\right\rangle \left\langle
e_{i},y\right\rangle +\sum_{j\in J}\left\langle x,f_{j}\right\rangle
\left\langle f_{j},y\right\rangle \right. -\left. 2\sum_{i\in I,j\in
J}\left\langle x,e_{i}\right\rangle \left\langle f_{j},y\right\rangle
\left\langle e_{i},f_{j}\right\rangle \right\vert  \label{2.10} \\
& \leq \frac{1}{2}\left\vert \left\langle x,y\right\rangle \right\vert
+\left\vert \sum_{i\in I}\left\langle x,e_{i}\right\rangle \left\langle
e_{i},y\right\rangle +\sum_{j\in J}\left\langle x,f_{j}\right\rangle
\left\langle f_{j},y\right\rangle \right.  \notag \\
& \qquad \qquad \qquad -\left. 2\sum_{i\in I,j\in J}\left\langle
x,e_{i}\right\rangle \left\langle f_{j},y\right\rangle \left\langle
e_{i},f_{j}\right\rangle -\frac{1}{2}\left\vert \left\langle
x,y\right\rangle \right\vert \right\vert  \notag \\
& \leq \frac{1}{2}\left[ \left\vert \left\langle x,y\right\rangle
\right\vert +\left\Vert x\right\Vert \left\Vert y\right\Vert \right] . 
\notag
\end{align}
\end{corollary}

\begin{proof}
The first inequality follows by the triangle inequality for the modulus. The
second inequality follows by (\ref{2.1}) on adding the quantity $\frac{1}{2}%
\left\vert \left\langle x,y\right\rangle \right\vert $ on both sides.
\end{proof}

\begin{remark}

\begin{enumerate}
\item \label{r2.4}If we choose in (\ref{2.10}), $x=y,$ then we get:%
\begin{align}
& \left\vert \sum_{i\in I}\left\vert \left\langle x,e_{i}\right\rangle
\right\vert ^{2}+\sum_{j\in J}\left\vert \left\langle x,f_{j}\right\rangle
\right\vert ^{2}-2\sum_{i\in I,j\in J}\left\langle x,e_{i}\right\rangle
\left\langle f_{j},x\right\rangle \left\langle e_{i},f_{j}\right\rangle
\right\vert  \label{2.11} \\
& \leq \left\vert \sum_{i\in I}\left\vert \left\langle x,e_{i}\right\rangle
\right\vert ^{2}+\sum_{j\in J}\left\vert \left\langle x,f_{j}\right\rangle
\right\vert ^{2}\right.  \notag \\
& \qquad \qquad \qquad -\left. 2\sum_{i\in I,j\in J}\left\langle
x,e_{i}\right\rangle \left\langle f_{j},x\right\rangle \left\langle
e_{i},f_{j}\right\rangle -\frac{1}{2}\left\Vert x\right\Vert ^{2}\right\vert
+\frac{1}{2}\left\Vert x\right\Vert ^{2}  \notag \\
& \leq \left\Vert x\right\Vert ^{2}.  \notag
\end{align}%
We observe that (\ref{2.11}) will generate Bessel's inequality if $\left\{
e_{i}\right\} _{i\in I},$ $\left\{ f_{j}\right\} _{j\in J}$ are disjoint
parts of a larger orthonormal family.

\item From (\ref{2.8}) one can obtain:%
\begin{equation}
\left\vert \frac{\left\langle x,e\right\rangle \left\langle e,y\right\rangle 
}{\left\Vert e\right\Vert ^{2}}+\frac{\left\langle x,f\right\rangle
\left\langle f,y\right\rangle }{\left\Vert f\right\Vert ^{2}}-2\cdot \frac{%
\left\langle x,e\right\rangle \left\langle f,y\right\rangle \left\langle
e,f\right\rangle }{\left\Vert e\right\Vert ^{2}\left\Vert f\right\Vert ^{2}}%
\right\vert \leq \frac{1}{2}\left[ \left\Vert x\right\Vert \left\Vert
y\right\Vert +\left\vert \left\langle x,y\right\rangle \right\vert \right]
\label{2.12}
\end{equation}%
and in particular 
\begin{equation}
\left\vert \frac{\left\vert \left\langle x,e\right\rangle \right\vert ^{2}}{%
\left\Vert e\right\Vert ^{2}}+\frac{\left\vert \left\langle x,f\right\rangle
\right\vert ^{2}}{\left\Vert f\right\Vert ^{2}}-2\cdot \frac{\left\langle
x,e\right\rangle \left\langle f,e\right\rangle \left\langle e,f\right\rangle 
}{\left\Vert e\right\Vert ^{2}\left\Vert f\right\Vert ^{2}}\right\vert \leq
\left\Vert x\right\Vert ^{2},  \label{2.13}
\end{equation}%
for any $x,y\in H.$
\end{enumerate}
\end{remark}

The case of real inner products will provide a natural genearlization for
Precupanu's inequality (\ref{1.1}):

\begin{corollary}
\label{c2.5}Let $\left( H;\left\langle \cdot ,\cdot \right\rangle \right) $
be a real inner product space and $\left\{ e_{i}\right\} _{i\in I},$ $%
\left\{ f_{j}\right\} _{j\in J}$ two finite families of orthonormal vectors
in $\left( H;\left\langle \cdot ,\cdot \right\rangle \right) .$ For any $%
x,y\in H\backslash \left\{ 0\right\} $ one has the double inequality:%
\begin{multline}
\frac{1}{2}\left[ \left\vert \left\langle x,y\right\rangle \right\vert
-\left\Vert x\right\Vert \left\Vert y\right\Vert \right]  \label{2.14} \\
\leq \sum_{i\in I}\left\langle x,e_{i}\right\rangle \left\langle
y,e_{i}\right\rangle +\sum_{j\in J}\left\langle x,f_{j}\right\rangle
\left\langle y,f_{j}\right\rangle -2\sum_{i\in I,j\in J}\left\langle
x,e_{i}\right\rangle \left\langle y,f_{j}\right\rangle \left\langle
e_{i},f_{j}\right\rangle \\
\leq \frac{1}{2}\left[ \left\Vert x\right\Vert \left\Vert y\right\Vert
+\left\vert \left\langle x,y\right\rangle \right\vert \right] .
\end{multline}%
In particular, we have%
\begin{align}
0& \leq \sum_{i\in I}\left\langle x,e_{i}\right\rangle ^{2}+\sum_{j\in
J}\left\langle x,f_{j}\right\rangle ^{2}-2\sum_{i\in I,j\in J}\left\langle
x,e_{i}\right\rangle \left\langle x,f_{j}\right\rangle \left\langle
e_{i},f_{j}\right\rangle  \label{2.15} \\
& \leq \left\Vert x\right\Vert ^{2},  \notag
\end{align}%
for any $x\in H.$
\end{corollary}

\begin{remark}
Similar particular inequalities to those incorporated in (\ref{2.7}) -- (\ref%
{2.13}) may be stated, but we omit them.
\end{remark}

\section{Refinements of Kurepa's Inequality\label{s3}}

Let $\left( H;\left\langle \cdot ,\cdot \right\rangle \right) $ be a real
inner product space generating the norm $\left\Vert \cdot \right\Vert .$ The 
\textit{complexification} $H_{\mathbb{C}}$ of $H$ is defined as a complex
linear space $H\times H$ of all ordered pairs $\left( x,y\right) ,$ $x,y\in
H $ endowed with the operations:%
\begin{align*}
\left( x,y\right) +\left( x^{\prime },y^{\prime }\right) & :=\left(
x+x^{\prime },y+y^{\prime }\right) ,\qquad x,x^{\prime },y,y^{\prime }\in H;
\\
\left( \sigma +i\tau \right) \cdot \left( x,y\right) & :=\left( \sigma
x-\tau y,\tau x+\sigma y\right) ,\qquad x,y\in H\text{ \ and \ }\sigma ,\tau
\in \mathbb{R}.
\end{align*}%
On $H_{\mathbb{C}}:=H\times H,$ endowed with the above operations, one can
now canonically define the \textit{scalar product} $\left\langle \cdot
,\cdot \right\rangle _{\mathbb{C}}$ by:%
\begin{equation}
\left\langle z,z^{\prime }\right\rangle _{\mathbb{C}}:=\left\langle
x,x^{\prime }\right\rangle +\left\langle y,y^{\prime }\right\rangle +i\left[
\left\langle x^{\prime },y\right\rangle -\left\langle x,y^{\prime
}\right\rangle \right]  \label{3.1}
\end{equation}%
where $z=\left( x,y\right) ,$ $z^{\prime }=\left( x^{\prime },y^{\prime
}\right) \in H_{\mathbb{C}}.$ Obviously,%
\begin{equation*}
\left\Vert z\right\Vert _{\mathbb{C}}^{2}=\left\Vert x\right\Vert
^{2}+\left\Vert y\right\Vert ^{2},\qquad z=\left( x,y\right) \in H_{\mathbb{C%
}}.
\end{equation*}%
One can also define the \textit{conjugate} of a vector $z=\left( x,y\right) $
by $\bar{z}:=\left( x,-y\right) .$ It is easy to see that, the elements of $%
H_{\mathbb{C}},$ under defined operations, behave as formal
\textquotedblleft complex\textquotedblright\ combinations $x+iy$ with $%
x,y\in H.$ Because of this, we may write $z=x+iy$ instead of $z=\left(
x,y\right) .$ Thus, $\bar{z}=x-iy.$

Under this setting, S. Kurepa \cite{K} proved the following refinement of
Schwarz's inequality:%
\begin{equation}
\left\vert \left\langle a,z\right\rangle _{\mathbb{C}}\right\vert ^{2}\leq 
\frac{1}{2}\left\Vert a\right\Vert ^{2}\left[ \left\Vert z\right\Vert _{%
\mathbb{C}}^{2}+\left\vert \left\langle z,\bar{z}\right\rangle _{\mathbb{C}%
}\right\vert \right] \leq \left\Vert a\right\Vert ^{2}\left\Vert
z\right\Vert _{\mathbb{C}}^{2},  \label{3.2}
\end{equation}%
for any $a\in H$ and $z\in H_{\mathbb{C}}.$

This was motivated by generalising the de Bruijn result for sequences of
real and complex numbers obtained in \cite{BR}.

The following result may be stated.

\begin{theorem}
\label{t3.1}Let $\left( H;\left\langle \cdot ,\cdot \right\rangle \right) $
be a real inner product space and $\left\{ e_{i}\right\} _{i\in I},\left\{
f_{j}\right\} _{j\in J}$ two finite families in $H.$ If $\left( H_{\mathbb{C}%
};\left\langle \cdot ,\cdot \right\rangle _{\mathbb{C}}\right) $ is the
complexification of $\left( H;\left\langle \cdot ,\cdot \right\rangle
\right) ,$ then for any $w\in H_{\mathbb{C}},$ we have the inequalities%
\begin{align}
& \left\vert \sum_{i\in I}\left\langle w,e_{i}\right\rangle _{\mathbb{C}%
}^{2}+\sum_{j\in J}\left\langle w,f_{j}\right\rangle _{\mathbb{C}%
}^{2}-2\sum_{i\in I,j\in J}\left\langle w,e_{i}\right\rangle _{\mathbb{C}%
}\left\langle w,f_{j}\right\rangle _{\mathbb{C}}\left\langle
e_{i},f_{j}\right\rangle \right\vert  \label{3.3} \\
& \leq \frac{1}{2}\left\vert \left\langle w,\bar{w}\right\rangle _{\mathbb{C}%
}\right\vert +\left\vert \sum_{i\in I}\left\langle w,e_{i}\right\rangle _{%
\mathbb{C}}^{2}+\sum_{j\in J}\left\langle w,f_{j}\right\rangle _{\mathbb{C}%
}^{2}\right.  \notag \\
& \qquad \qquad \qquad -\left. 2\sum_{i\in I,j\in J}\left\langle
w,e_{i}\right\rangle _{\mathbb{C}}\left\langle w,f_{j}\right\rangle _{%
\mathbb{C}}\left\langle e_{i},f_{j}\right\rangle -\frac{1}{2}\left\langle w,%
\bar{w}\right\rangle _{\mathbb{C}}\right\vert  \notag \\
& \leq \frac{1}{2}\left[ \left\Vert w\right\Vert _{\mathbb{C}%
}^{2}+\left\vert \left\langle w,\bar{w}\right\rangle _{\mathbb{C}%
}\right\vert \right] \leq \left\Vert w\right\Vert _{\mathbb{C}}^{2}.  \notag
\end{align}
\end{theorem}

\begin{proof}
Define $g_{j}\in H_{\mathbb{C}},$ $g_{j}:=\left( e_{j},0\right) ,$ $j\in I.$
For any $k,j\in I$ we have%
\begin{equation*}
\left\langle g_{k},g_{j}\right\rangle _{\mathbb{C}}=\left\langle \left(
e_{k},0\right) ,\left( e_{j},0\right) \right\rangle _{\mathbb{C}%
}=\left\langle e_{k},e_{j}\right\rangle =\delta _{kj},
\end{equation*}%
therefore $\left\{ g_{j}\right\} _{j\in I}$ is an orthonormal family in $%
\left( H_{\mathbb{C}};\left\langle \cdot ,\cdot \right\rangle _{\mathbb{C}%
}\right) .$

If we apply Corollary \ref{c2.3} for $\left( H_{\mathbb{C}};\left\langle
\cdot ,\cdot \right\rangle _{\mathbb{C}}\right) ,$ $x=w,$ $y=\bar{w},$ we
may write:%
\begin{align}
& \left\vert \sum_{i\in I}\left\langle w,e_{i}\right\rangle _{\mathbb{C}%
}\left\langle e_{i},\bar{w}\right\rangle _{\mathbb{C}}+\sum_{j\in
J}\left\langle w,f_{j}\right\rangle \left\langle f_{j},\bar{w}\right\rangle
\right.  \label{3.4} \\
& \qquad \qquad \qquad -\left. 2\sum_{i\in I,j\in J}\left\langle
w,e_{i}\right\rangle _{\mathbb{C}}\left\langle f_{j},\overline{w}%
\right\rangle _{\mathbb{C}}\left\langle e_{i},f_{j}\right\rangle \right\vert
\notag \\
& \leq \frac{1}{2}\left\Vert w\right\Vert _{\mathbb{C}}\left\Vert \bar{w}%
\right\Vert _{\mathbb{C}}+\left\vert \sum_{i\in I}\left\langle
w,e_{i}\right\rangle _{\mathbb{C}}\left\langle e_{i},\bar{w}\right\rangle _{%
\mathbb{C}}+\sum_{j\in J}\left\langle w,f_{j}\right\rangle \left\langle
f_{j},\bar{w}\right\rangle \right.  \notag \\
& \qquad \qquad \qquad -\left. 2\sum_{i\in I,j\in J}\left\langle
w,e_{i}\right\rangle _{\mathbb{C}}\left\langle f_{j},\overline{w}%
\right\rangle _{\mathbb{C}}\left\langle e_{i},f_{j}\right\rangle -\frac{1}{2}%
\left\langle w,\bar{w}\right\rangle _{\mathbb{C}}\right\vert  \notag \\
& \leq \frac{1}{2}\left[ \left\vert \left\langle w,\bar{w}\right\rangle _{%
\mathbb{C}}\right\vert +\left\Vert w\right\Vert _{\mathbb{C}}\left\Vert \bar{%
w}\right\Vert _{\mathbb{C}}\right] .  \notag
\end{align}%
However, for $w:=\left( x,y\right) \in H_{\mathbb{C}}$, we have $\bar{w}%
=\left( x,-y\right) $ and%
\begin{equation*}
\left\langle e_{j},\bar{w}\right\rangle _{\mathbb{C}}=\left\langle \left(
e_{j},0\right) ,\left( x,-y\right) \right\rangle _{\mathbb{C}}=\left\langle
e_{j},x\right\rangle +i\left\langle e_{j},y\right\rangle
\end{equation*}%
and%
\begin{equation*}
\left\langle w,e_{j}\right\rangle _{\mathbb{C}}=\left\langle \left(
x,y\right) ,\left( e_{j},0\right) \right\rangle _{\mathbb{C}}=\left\langle
x,e_{j}\right\rangle +i\left\langle e_{j},y\right\rangle
\end{equation*}%
showing that $\left\langle e_{j},\bar{w}\right\rangle =\left\langle
w,e_{j}\right\rangle $ for any $j\in I$. A similar relation is true for $%
f_{j}$ and since%
\begin{equation*}
\left\Vert w\right\Vert _{\mathbb{C}}=\left\Vert \bar{w}\right\Vert _{%
\mathbb{C}}=\left( \left\Vert x\right\Vert ^{2}+\left\Vert y\right\Vert
^{2}\right) ^{\frac{1}{2}},
\end{equation*}%
hence from (\ref{3.4}) we deduce the desired inequality (\ref{3.3}).
\end{proof}

\begin{remark}
It is obvious that, if one family, say $\left\{ f_{j}\right\} _{j\in J}$ is
empty, then, on observing that all sums $\sum_{j\in J}$ should be zero, from
(\ref{3.3}) one would get \cite{SSD}%
\begin{align}
& \left\vert \sum_{i\in I}\left\langle w,e_{i}\right\rangle _{\mathbb{C}%
}^{2}\right\vert  \label{3.5} \\
& \leq \frac{1}{2}\left\vert \left\langle w,\bar{w}\right\rangle _{\mathbb{C}%
}\right\vert +\left\vert \sum_{i\in I}\left\langle w,e_{i}\right\rangle _{%
\mathbb{C}}^{2}-\frac{1}{2}\left\langle w,\bar{w}\right\rangle _{\mathbb{C}%
}\right\vert  \notag \\
& \leq \frac{1}{2}\left[ \left\Vert w\right\Vert _{\mathbb{C}%
}^{2}+\left\vert \left\langle w,\bar{w}\right\rangle _{\mathbb{C}%
}\right\vert \right] \leq \left\Vert w\right\Vert _{\mathbb{C}}^{2}.  \notag
\end{align}%
If in (\ref{3.5}) one assumes that the family $\left\{ e_{i}\right\} _{i\in
I}$ contains only one element $e=\frac{a}{\left\Vert a\right\Vert },a\neq 0,$
then by selecting $w=z,$ one would deduce 
\begin{eqnarray*}
\left\vert \left\langle a,z\right\rangle _{\mathbb{C}}\right\vert ^{2} &\leq
&\left\vert \left\langle a,z\right\rangle _{\mathbb{C}}^{2}-\frac{1}{2}%
\left\langle z,\bar{z}\right\rangle _{\mathbb{C}}\right\vert +\frac{1}{2}%
\left\vert \left\langle z,\bar{z}\right\rangle _{\mathbb{C}}\right\vert \\
&\leq &\frac{1}{2}\left\Vert a\right\Vert ^{2}\left[ \left\Vert z\right\Vert
_{\mathbb{C}}^{2}+\left\vert \left\langle z,\bar{z}\right\rangle _{\mathbb{C}%
}\right\vert \right] ,
\end{eqnarray*}%
which is a refinement for Kurepa's inequality (\ref{3.2}).
\end{remark}

\begin{acknowledgement}
The author would like to thank the anonymous referee for his/her comments
that have been implemented in the final version of this paper.
\end{acknowledgement}

\end{document}